\documentclass[graybox]{svmult}

\usepackage{type1cm}        
\usepackage{makeidx}         
\usepackage{graphicx}        
\usepackage{multicol}        
\usepackage[bottom]{footmisc}

\usepackage{newtxtext}       %
\usepackage[varvw]{newtxmath}       

\usepackage{tikz-cd}

\makeatletter
\newcommand{\ostar}{\mathbin{\mathpalette\make@circled\ast}}
\newcommand{\make@circled}[2]{%
  \ooalign{$\m@th#1\smallbigcirc{#1}$\cr\hidewidth$\m@th#1#2$\hidewidth\cr}%
}
\newcommand{\smallbigcirc}[1]{%
  \vcenter{\hbox{\scalebox{0.77778}{$\m@th#1\bigcirc$}}}%
}
\makeatother

\makeindex             

\begin{document}

\title*{Intersection Multiplicity in Loop Spaces and the String Topology Coproduct}
\author{Philippe Kupper 
 and\\ Maximilian Stegemeyer 
 }
\institute{Philippe Kupper \at Institut für Algebra und Geometrie, {Karlsruher Institut für Technologie}, {Karlsruhe}, {76128},  Germany, \email{philippe.kupper@kit.edu}
\and Maximilian Stegemeyer \at Mathematisches Institut, Universität Freiburg, Ernst-Zermelo-Straße 1, Freiburg, 79104, Germany, \email{maximilian.stegemeyer@math.uni-freiburg.de}}
%
%
\maketitle

\abstract*{The string topology coproduct is often perceived as a counterpart in string topology to the Chas-Sullivan product. However, in certain aspects the string topology coproduct is much harder to understand than the Chas-Sullivan product. In particular the coproduct is not homotopy-invariant and it seems much harder to compute.
In this article we give an overview over the string topology coproduct and use the notion of intersection multiplicity of homology classes in loop spaces to show that the string topology coproduct and the based string topology coproduct are trivial for certain classes of manifolds. 
In particular we show that the string topology coproduct vanishes on product manifolds where both factors have vanishing Euler characteristic and we show that the based coproduct is trivial for total spaces of fiber bundles with sections.
We also discuss implications of these results.}

\abstract{The string topology coproduct is often perceived as a counterpart in string topology to the Chas-Sullivan product. However, in certain aspects the string topology coproduct is much harder to understand than the Chas-Sullivan product. In particular the coproduct is not homotopy-invariant and it seems much harder to compute.
In this article we give an overview over the string topology coproduct and use the notion of intersection multiplicity of homology classes in loop spaces to show that the string topology coproduct and the based string topology coproduct are trivial for certain classes of manifolds. 
In particular we show that the string topology coproduct vanishes on product manifolds where both factors have vanishing Euler characteristic and we show that the based coproduct is trivial for total spaces of fiber bundles with sections.
We also discuss implications of these results.}

\section{Introduction}\label{sec1}

String topology is the study of algebraic structures on the homology of the free loop space of a closed manifold $M$.
The most prominent operation is the \textit{Chas-Sullivan product} 
$$  \wedge\colon \mathrm{H}_i(\Lambda M)\otimes \mathrm{H}_j(\Lambda M)\to \mathrm{H}_{i+j-n}(\Lambda M).    $$
Here, $\Lambda M$ is the free loop space of $M $ and $n$ is the dimension of $M$.
This product was first introduced by Chas and Sullivan in \cite{chas:1999}.
Its geometric idea is to concatenate loops in $M$ that share the same basepoint.

The Chas-Sullivan product is quite well-understood.
Indeed, it is known that it is a homotopy invariant of the underlying manifold $M$, see \cite{cohen:2008} and there are many computations of the Chas-Sullivan product, see e.g. \cite{cohen:2004}, \cite{hepworth:2009}, \cite{hepworth:2010}, \cite{kupper:2021} and \cite{tamanoi:2006} just to name a few.

Following ideas by Sullivan in \cite{sullivan:2004}, Goresky and Hingston \cite{goresky:2009} define the \textit{string topology coproduct}
$$   \vee\colon \mathrm{H}_i(\Lambda M,M) \to ( \mathrm{H}_{\bullet}(\Lambda M,M)\otimes \mathrm{H}_{\bullet}(\Lambda M,M)  )_{i+1-n}   .  $$
Here, we consider $M$ as a subspace of $\Lambda M$ via the embedding of $M$ as the constant loops.
It turns out that the string topology coproduct is harder to understand than the Chas-Sullivan product.
First of all, the string topology coproduct is not a homotopy invariant of the underlying manifold.
Naef shows in \cite{naef:2021} that the string topology coproduct can tell apart two Lens spaces which are homotopy equivalent, but not homeomorphic.
It is also very hard to understand what structure the Chas-Sullivan product and the string topology coproduct form together.
For results in this direction we refer to \cite{cieliebak2020poincar}, \cite{naef2019string} and \cite{rivera2019singular}.
Moreover, it seems that the string topology coproduct is harder to compute than the Chas-Sullivan product.
There are computations for spheres and projective spaces, see \cite{goresky:2009}, \cite{hingston:2017} and \cite{stegemeyer:2022}.
The second author has shown in \cite{stegemeyer:2021} that the string topology coproduct is trivial for all compact simply connected Lie groups of rank at least $2$.
Naito shows in \cite{naef:2021} that the string topology coproduct with rational coefficients is trivial for large classes of spaces.
This raises the question for which class of spaces the string topology coproduct is actually non-trivial.

The purpose of this article is two-fold.
In Section \ref{sec_survey} we give an overview over the string topology coproduct, its failure to be a homotopy invariant and in particular about computational examples as well as relations to Riemannian geometry.
We shall see that the list of examples where the string topology coproduct has been computed is not very long and that there are only very few examples where the string topology coproduct is non-trivial.
In Section \ref{sec_int_mult} and \ref{sec_results} we shall then see how the concept of \textit{intersection multiplicity} can be used to deduce further triviality statements for the string topology coproduct and for the based version of the string topology coproduct.
While the proofs of our results in Section \ref{sec_results} are quite elementary we found it interesting that one can already use very naive arguments and the notion of intersection multiplicity to show that the string topology coproduct and its based counterpart vanish in certain situations.

In particular we shall show that the string topology coproduct vanishes for certain product manifolds.
\begin{unnumberdtheorem}[Theorem \ref{theorem_int_mult_products}]
    Let $R$ be a commutative unital ring and let $M = N\times P$ be a product of two closed $R$-oriented manifolds each with vanishing Euler characteristic.
    Take homology with coefficients in $R$. Then the string topology coproduct on $M$ is trivial.
\end{unnumberdtheorem}
Moreover, we study the based string topology coproduct on the total space of a fiber bundle with a section.
\begin{unnumberdtheorem}[Theorem \ref{tristan} and Corollary \ref{isolde}]
Let $R$ be a commutative unital ring and take homology with coefficients in $R$.   
    Let $p\colon E\to B$ be a smooth fiber bundle with a section, where $B$ and $E$ are closed manifolds. Assume that both $B$ and the fiber are positive dimensional. Assume further that $E$ is $R$-oriented.
    Then the based string topology coproduct on $E$ is trivial.
    
    Conversely, if a closed $R$-oriented manifold $M$ has non-trivial based string topology coproduct then it cannot be the total space of a fiber bundle with a section, where both base and fiber are positive dimensional.
\end{unnumberdtheorem}
Finally, we shall see that a non-trivial string topology coproduct can be used to infer some properties of maps in homology.
\begin{unnumberdtheorem}[Corollary \ref{cor_submanifold_coproduct} and Theorem \ref{thm_cpn}]
    Let $R$ be a commutative unital ring and take homology with coefficients in $R$. 
    Let $M$ be an $R$-oriented closed manifold and $N$ a closed submanifold of $M$ such that the normal bundle of $N$ has a nowhere vanishing section.
    \begin{enumerate}
    \item The string topology coproduct vanishes on the subspace $U \subseteq \mathrm{H}_{\bullet}(\Lambda M,M)$ which is the image of the map
    $$  i_N\colon  \mathrm{H}_{\bullet}(\Lambda N) \xrightarrow[]{j_*} \mathrm{H}_{\bullet}(\Lambda M)\to \mathrm{H}_{\bullet}(\Lambda M,M) $$
     induced by the inclusion $j\colon \Lambda N\hookrightarrow\Lambda M$.
    \item For $M =\mathbb{C}P^n$ and $N = \mathbb{C}P^m$ with $m< \tfrac{n}{2}$ the map
    $$   j_* \colon \mathrm{H}_k(\Lambda \mathbb{C}P^m;\mathbb{Q}) \to \mathrm{H}_k(\Lambda \mathbb{C}P^n;\mathbb{Q})    $$
    in rational homology is trivial for $k\geq 4n+1$, where $j\colon \Lambda \mathbb{C}P^m\hookrightarrow \Lambda \mathbb{C}P^n$ is the inclusion.
    \end{enumerate}
\end{unnumberdtheorem}
\medskip


\noindent \textbf{Acknowledgements:} The authors would like to thank Florian Kranhold for a very interesting and helpful discussion about principal fibrations and confirming that no assumption on the fundamental groups is needed in Lemma \ref{lemma_principal_fibration}.

M. S. was partially funded by the Deutsche Forschungsgemeinschaft (German Research Foundation) -- grant agreement number 518920559.
M.S. is grateful for the support by the Danish National Research Foundation through the Copenhagen Centre for Geometry and Topology (DNRF151).

\section{The String Topology Coproduct}\label{sec_survey}

In this section we give an overview over the string topology coproduct.
We shall begin by reviewing its definition as well as the definitions of some related structrures like the Goresky-Hingston product and the based string topology coproduct.
We then discuss the algebraic properties and the question of homotopy invariance of the string topology coproduct.
Finally, we shall see some examples as well as connections to geometric questions about closed geodesics.

This exposition must stay brief in certain places but we will of course give as many references as possible to the existing literature.

\subsection{Definition and algebraic properties}\label{subsec_definition}

We begin by introducing the string topology coproduct.
Before we come to the particular definition we need to set the ground by introducing the free loop space of a closed manifold.

Let us consider a closed manifold $M$ and denote the unit interval by $I = [0,1]$.
We define the \textit{free loop space of $M$} to be $$  \Lambda M = \big\{ \gamma : I\to M\,|\, \gamma(0) = \gamma(1),\,\, \gamma \,\,\text{absolutely continuous}, \,\, \int_0^1 |\Dot{\gamma}(t)|^2 \,\mathrm{d}t < \infty \big\} \,   .$$
Here the norm of the tangent vector $|\Dot{\gamma}(t)|$, $t\in I$, is taken with respect to an arbitrary Riemannian metric, since if the integral $\int_0^1 |\Dot{\gamma}(t)|^2 \,\mathrm{d}t$ is finite for one Riemannian metric then it is finite for any Riemannian metric.
For the notion of absolute continuity of curves in a smooth manifold we refer to \cite[Definition 2.3.1]{klingenberg:1995}.
It turns out that the free loop space $\Lambda M$ can be given the structure of a Hilbert manifold, see \cite[Section 2.3]{klingenberg:1995}.
The base manifold $M$ is a submanifold of $\Lambda M$ which is embedded in $\Lambda M$ as the trivial loops, see \cite[Proposition 1.4.6]{klingenberg:78}.
If we choose a basepoint $p_0\in M$ we define
$$ \Omega_{p_0}M = \{ \gamma\in \Lambda M\,|\, \gamma(0) = p_0\}     $$
and call this space the \textit{based loop space of} $M$ \textit{at the point} $p_0$.
The based loop space is a submanifold of the free loop space $\Lambda M$.

For technical reasons it is very convenient to work with absolutely continuous loops.
However, we remark that the free loop space as we defined it above is homotopy equivalent to the continuous mapping space $C^0(\mathbb{S}^1,M)$ equipped with the compact-open topology as well as to the Fréchet manifold of piecewise smooth loops.
This was proven by Palais, see \cite{palais:1968}.
See also \cite[Section 2]{chataur:2015} for an overview about the different versions of the free loop space.

We note that the evaluation $\mathrm{ev}\colon \Lambda M\to M$ given by $\mathrm{ev}(\gamma) = \gamma(0)$ is a fibration with typical fiber $\mathrm{ev}^{-1}(\{p_0\}) = \Omega_{p_0} M$.
We call this the \textit{free loop fibration of} $M$.
For more properties of the free loop space in general and in particular of the free loop fibration we again refer to \cite{chataur:2015}.

We now turn to string topology.
The string topology coproduct goes back to ideas by Sullivan in \cite{sullivan:2004}.
The geometric idea is to cut apart loops which have self-intersections.
More precisely, for $s\in I$ consider the spaces 
$$    F_s = \{ \gamma\in \Lambda M\,|\, \gamma(0) = \gamma(s)\} \quad \text{and} \quad F = \{ (\gamma,s)\in \Lambda M \times I \,|\, \gamma(0) = \gamma(s) \}   .   $$
We have an obvious inclusion $F_s \hookrightarrow F$ as well as a map $\mathrm{cut}\colon F\to \Lambda M\times \Lambda M$ which takes $(\gamma,s)$ to the restrictions $(\gamma|_{[0,s]},\gamma|_{[s,1]})$.
The restrictions are appropriately re-parametrized so that they are defined on the unit interval.
This map restricts to a map $\mathrm{cut}_s\colon F_s \to \Lambda M\times \Lambda M$.
Hence, we obtain diagrams
\begin{equation} \label{eq_diagram_1}
    \begin{tikzcd}
    \Lambda M   \arrow[r, hookleftarrow]    & F_s     \arrow[]{r}{\mathrm{cut }} & \Lambda M \times \Lambda M
\end{tikzcd}
\end{equation}
as well as
\begin{equation} \label{eq_diagram_2}
\begin{tikzcd}
    \Lambda M \times I  \arrow[r, hookleftarrow]    & F     \arrow[]{r}{\mathrm{cut }} & \Lambda M \times \Lambda M  .
\end{tikzcd}
\end{equation}
In order to get interesting operations in homology we need to find morphisms in homology which - in case of the first diagram - map the homology $\mathrm{H}_{\bullet}(\Lambda M)$ to $\mathrm{H}_{\bullet}(F_s)$.
More generally, one is interested in the following situation. Let $R$ be a commutative unital ring and take homology with coefficients in $R$. Let $X$ be a Hilbert manifold and $Y\subseteq X$ a submanifold of finite codimension with $R$-orientable normal bundle.
Then one wants to introduce a morphism $\mathrm{H}_i(X)\to\mathrm{H}_{i-k}(Y)$ where $k$ is the codimension of $Y$ in $X$.
This can in fact be achieved by \textit{Gysin morphisms}, see e.g. \cite[Appendix B]{goresky:2009}.
If the manifold $M$ is $R$-oriented and if we take the first diagram then we can use a Gysin morphism which maps $\mathrm{H}_i(\Lambda M)\to \mathrm{H}_{i-n}(F_s)$. We can compose this with the cutting map to obtain a coproduct
$$    \vee_s \colon \mathrm{H}_i(\Lambda M)\to \mathrm{H}_{i-n}(\Lambda M\times \Lambda M) .    $$
However, it turns out that this coproduct is trivial except possibly in low degrees.
This was first observed by Tamanoi \cite{tamanoi2010loop}, see also  \cite{cohen2004polarized}, \cite{godin2007higher} and \cite[Section 9.1]{goresky:2009}.
Goresky and Hingston therefore use the diagram \eqref{eq_diagram_2} to define a coproduct as follows.
Let $[I]\in \mathrm{H}_1(I,\partial I)$ be a positively oriented generator.
Taking the homology cross product with $[I]$ defines a map $$ f_1 \colon \mathrm{H}_i(\Lambda M,M)\to \mathrm{H}_{i+1}(\Lambda M\times I, \Lambda M\times \partial I\cup M\times I).$$
Then one constructs a map 
$$f_2\colon \mathrm{H}_{i+1}(\Lambda M\times I, \Lambda M\times \partial I\cup M\times I)\to \mathrm{H}_{i+1-n}(F,\Lambda M\times \partial I\cup M\times I) . $$ 
Note that the map $f_2$ is defined for relative homology and the space $F$ is not an honest submanifold of $\Lambda M\times I$.
Therefore the map $f_2$ is more intricate than the Gysin morphisms that we saw above.
For details of the construction of the morphism $f_2$ we refer to \cite[Sections 8 and 9]{goresky:2009} and to \cite[Sections 1 and 2]{hingston:2017}.
Finally, note that the cutting map $\mathrm{cut}\colon F\to \Lambda M\times \Lambda M$ induces a map of pairs
$$    \mathrm{cut}\colon (F, \Lambda M\times \partial I\cup M\times I) \to (\Lambda M\times \Lambda M,\Lambda M\times M\cup M\times \Lambda M) .       $$
\begin{definition}
    Let $R$ be a commutative unital ring and $M$ an $R$-oriented closed manifold. Take homology with coefficients in $R$. The \textit{unsigned string topology coproduct} is defined as the composition 
$$     \vee = \mathrm{cut}_*\circ f_2\circ f_1\colon \mathrm{H}_i(\Lambda M,M)\to \mathrm{H}_{i+1-n}(\Lambda M\times \Lambda M,\Lambda M\times M\cup M\times \Lambda M) .       $$
\end{definition}
If we take homology with coefficients in a field we can compose this with the Künneth morphism to obtain a map
$$      \vee\colon \mathrm{H}_{i}(\Lambda M,M)   \to \big(   \mathrm{H}_{\bullet}(\Lambda M,M)\otimes \mathrm{H}_{\bullet}(\Lambda M,M)    \big)_{i+1-n}     .  $$

In the following we will see that in order to guarantee good algebraic properties one needs to introduce additional signs.
Therefore we call the above map the \textit{unsigned} string topology coproduct.

Similarly to the string topology coproduct one can define a \textit{based string topology coproduct} on the homology of the based loop space $\Omega_{p_0} M$.
Indeed the construction is just a restriction of the string topology coproduct.
Instead of the space $F$ we take
$$    F \cap (\Omega_{p_0}M \times I) = \{ (\gamma,s)\in \Omega_{p_0} M\times I\,|\, \gamma(s) = p_0\} .     $$
This gives a map, the \textit{based string topology coproduct}
$$    \vee_{\Omega}\colon \mathrm{H}_i(\Omega M,\{p_0\}) \to \mathrm{H}_{i+1-n}(\Omega M\times \Omega M, \Omega M\times \{p_0\}\cup \{p_0\}  \times \Omega M) .    $$
As one expects one gets the following compatibility.
The diagram
$$
\begin{tikzcd}    \mathrm{H}_i(\Omega ,\{p_0\}) \arrow[]{r}{j_*} \arrow[]{d}{\vee_{\Omega}} & \mathrm{H}_i(\Lambda ,M) \arrow[]{d}{\vee}
    \\
    \mathrm{H}_{i+1-n}(\Omega \times \Omega,\Omega\times \{ p_0\}\cup  \{p_0\}\times \Omega) \arrow[]{r}{\hphantom{l}(j\times j)_*\hphantom{bl}} & \mathrm{H}_{i+1-n}(\Lambda \times \Lambda,\Lambda\times M\cup M\times \Lambda) 
\end{tikzcd}
$$
commutes, see e.g. \cite[Proposition 2.4]{stegemeyer:2021}. 
Here, we use the notation $\Omega = \Omega_{p_0}M $ and $\Lambda = \Lambda M$ and $j\colon \Omega\to \Lambda $ is the inclusion.

If we take coefficients in a field $\mathbb{K}$ and if the homology of $\Lambda M$ is of finite type we directly obtain a dual product in cohomology.
More precisely, let $\alpha\in \mathrm{H}^i(\Lambda M,M;\mathbb{K})$ and $\beta\in \mathrm{H}^j(\Lambda  M,M ;\mathbb{K})$ be relative cohomology classes.
The \textit{unsigned Goresky-Hingston product} $\alpha\ostar \beta$ is defined to be the unique cohomology class in $\mathrm{H}^{i+j+n-1}(\Lambda M,M ;\mathbb{K})$ such that
$$   \langle \alpha \ostar \beta , X\rangle = \langle \alpha \times \beta , \vee X\rangle \quad \text{for all} \,\,\, X\in\mathrm{H}_{\bullet}(\Lambda M,M;\mathbb{K})\, .    $$
Here $\langle\cdot,\cdot\rangle$ is the Kronecker pairing between homology and cohomology which - by our assumption of $\Lambda M$ being of finite type - is a perfect pairing.
Similarly, one can define a \textit{based Goresky-Hingston product} on the cohomology of the based loop space.

The Goresky-Hingston product was first defined intrinsically on the cohomology of the free loop space in \cite{goresky:2009}.
Hingston and Wahl show in \cite{hingston:2017} that the product which we defined above as the dual to the string topology coproduct in fact agrees with this intrinsically defined Goresky-Hingston product.

We now state the algebraic properties of the string topology coproduct and of the Goresky-Hingston product.
In \cite{hingston:2017} Hingston and Wahl introduce signs in the string topology coproduct in the following way.
Taking singular chains they define a chain-level version of the unsigned coproduct $$\vee_{\mathrm{chain}}\colon \mathrm{C}_{\bullet}(\Lambda M,M) \to \mathrm{C}_{\bullet}(\Lambda M,M) \otimes \mathrm{C}_{\bullet}(\Lambda M,M) .  $$
Then they define a \textit{signed} version of the chain-level string topology coproduct by setting
$$     \vee_{\mathrm{sign}} X =    \sum (-1)^{n + np}  X^0_p \otimes X^1_q       $$
where $\vee_{\mathrm{chain}} X =\sum   X^0_p \otimes X^1_q $.
This coproduct $\vee_{\mathrm{sign}}$ then defines a homology coproduct 
$$  \vee_{\mathrm{sign}} \colon \mathrm{H}_i(\Lambda M,M) \to \mathrm{H}_{i+1-n}(\Lambda M\times \Lambda M,\Lambda M\times M\cup M\times \Lambda M)     $$
which we call the \textit{(signed) string topology coproduct}.
Moreover if we take homology with coefficients in a field and if we assume that the homology of $\Lambda M$ is of finite type then we can again take the dual product in cohomology
$$   \ostar_{\mathrm{sign}}  \colon \mathrm{H}^{i}(\Lambda M,M) \otimes \mathrm{H}^j(\Lambda M,M) \to \mathrm{H}^{i+j+n-1}(\Lambda M,M)    $$
which we call the \textit{signed Goresky-Hingston product}. 

With these signs it turns out that the string topology coproduct and the Goresky-Hingston product have good algebraic properties which will be summed up in the following theorem.
For a proof we refer to \cite[Section 2]{hingston:2017}.

\begin{theorem}[Theorem 2.14 in \cite{hingston:2017}]
    Let $M$ be a closed oriented manifold and take homology with coefficients in a field $\mathbb{K}$.
    \begin{enumerate}
        \item The string topology coproduct 
        $$  \vee_{\mathrm{sign}} \colon  \mathrm{H}_{i-n}(\Lambda M,M) \to   \bigoplus_{j + l = i +1}  \mathrm{H}_{j-n}(\Lambda M,M) \otimes \mathrm{H}_{l-n}(\Lambda M,M)    $$
        is the twisted suspension of a graded coassociative and cocommutative coproduct on the degree-shifted
        homology $\mathrm{H}_{\bullet-n}(\Lambda M,M)$.
        \item The Goresky-Hingston product 
        $$     \ostar_{\mathrm{sign}} \colon \mathrm{H}^{i-n}(\Lambda M,M)\otimes \mathrm{H}^{j-n}(\Lambda M,M) \to \mathrm{H}^{i+j-1-n}      (\Lambda M,M)  $$
        is the twisted desuspension of a graded associative and commutative product on the degree-shifted
    cohomology $\mathrm{H}^{\bullet -n} (\Lambda M,M)$.
    \end{enumerate}
\end{theorem}
For a definition of the above algebraic notions of twisted suspension and twisted desuspension we refer to \cite{hingston:2017} and \cite{klamt2013natural}.
\begin{remark}
    We note as in \cite[Theorem 2.14]{hingston:2017} that the Goresky-Hingston product in particular satisfies the following signed commutativity and associativity relations.
    For $$\varphi, \psi,\chi\in \mathrm{H}^{\bullet}(\Lambda M,M)    \quad \text{with} \quad |\varphi| = i, \,\,|\psi| = j,\,\, |\chi| = k  $$
    we have
$$       \varphi\ostar_{\mathrm{sign}} \psi =  (-1)^{(i+n)(j+n)+1} \psi \ostar_{\mathrm{sign}} \varphi $$ and $$(\varphi\ostar_{\mathrm{sign}}\psi)\ostar_{\mathrm{sign}} \chi = (-1)^{k+n+1} \varphi\ostar_{\mathrm{sign}} (\psi\ostar_{\mathrm{sign}}\chi) .       $$
Moreover, we note that the Goresky-Hingston product is not unital.
\end{remark}
Finally, we remark that in Sections \ref{sec_int_mult} and \ref{sec_results} we will always use the unsigned version of the string topology coproduct, because we will mostly deal with situations where the string topology coproduct is trivial.
Of course, the unsigned version of the coproduct is trivial if and only if the signed version is trivial.

\subsection{Homotopy Invariance}\label{subsec_invariance}

It is now an obvious question how the string topology coproduct behaves under homotopy equivalences of the underlying manifold.
More precisely, let $f\colon M\to N$ be a map between closed manifolds and take homology with coefficients in a commutative unital ring $R$.
Then there is an induced map of loop spaces $\Lambda f\colon \Lambda M\to \Lambda N$.
It is easy to see that if $f$ is a homotopy equivalence then $\Lambda f$ is a homotopy equivalence as well.
Consequently, as $R$-modules the homologies $\mathrm{H}_{\bullet}(\Lambda M)$ and $\mathrm{H}_{\bullet}(\Lambda N)$ are isomorphic.
It is an important question if the respective string topology operations are preserved under these isomorphisms.
First, we note that the Chas-Sullivan product indeed behaves well under homotopy equivalences.
\begin{theorem}[Theorem 1 in \cite{cohen:2008}]
    Let $R$ be a commutative unital ring and let $M$ and $N$ be closed $R$-oriented manifolds.
    Let $f\colon M\to N$ be an orientation-preserving homotopy equivalence.
    Then the induced map $$(\Lambda f)_*   \colon \mathrm{H}_{\bullet}(\Lambda M)\to \mathrm{H}_{\bullet}(\Lambda N)    $$
    in homology with coefficients in $R$ is an isomorphism of rings if we equip $\mathrm{H}_{\bullet}(\Lambda M)$ and $\mathrm{H}_{\bullet}(\Lambda N)$ with the respective Chas-Sullivan product.
\end{theorem}

It turns out that the situation is much more complicated when one considers the string topology coproduct.
Indeed, Naef shows in \cite{naef:2021} that there is a homotopy equivalence $f\colon L(7,1)\to L(7,2)$ between the three-dimensional lens spaces $L(7,1)$ and $L(7,2)$ which does not preserve the string topology coproduct.
Naef also argues that this non-homotopy invariance is related to Whitehead torsion, see \cite[Section 2.4]{naef:2021}.
See also \cite{naef:2022} for a discussion of this example.
However, under certain additional assumptions on either the involved manifolds or on the homotopy equivalences one obtains statements about homotopy invariance.
\begin{theorem}[Corollary 1.2 in \cite{rivera2022invariance}]
    Let $f\colon M\to N$ be a homotopy equivalence between simply connected closed manifolds.
    Then the map 
    $$(\Lambda f)^*   \colon \mathrm{H}^{\bullet}(\Lambda N, N ;\mathbb{Q})\to \mathrm{H}^{\bullet}(\Lambda M, M;\mathbb{Q})    $$
    is an isomorphism of the respective Goresky-Hingston algebras. 
\end{theorem}
See also \cite{naef2019string}.
\begin{theorem}[Theorem A in \cite{hingston:2019}]
    Let $M, N$ be closed oriented manifolds and $f\colon M\to N$ be a homotopy equivalence homotopic
    to a $\Delta$–transverse map satisfying either condition B or condition D, or to a map satisfying condition D'.
    Then the induced map
    $$(\Lambda f)_*   \colon \mathrm{H}_{\bullet}(\Lambda M,M;\mathbb{Z})\to \mathrm{H}_{\bullet}(\Lambda N,N;\mathbb{Z})    $$
    intertwines the string topology coproduct on $M$ and $N$, respectively, up to a sign.
\end{theorem}
For the notion of $\Delta$\textit{-transverse map} and the conditions B,D and D' we refer to \cite{hingston:2019}.
We end this section by noting that the based string topology coproduct is preserved under homotopy equivalences.
If $f\colon M\to N$ is a map between pointed spaces, then there is an induced map between the based loop spaces $\Omega f\colon \Omega M\to \Omega N$.
Nancy Hingston shows the following.
\begin{theorem}[Corollary in \cite{hingston:2010}]
    Let $f\colon M\to N$ be a based orientation-preserving homotopy equivalence between closed simply connected manifolds.
    Then the induced map
    $$(\Omega f)_*   \colon \mathrm{H}_{\bullet}(\Omega M,\{p_0\};\mathbb{Z})\to \mathrm{H}_{\bullet}(\Omega N,\{q_0\};\mathbb{Z})    $$
    intertwines the based string topology coproduct on $M$ and $N$.
\end{theorem}

\subsection{Examples and (Non-)Triviality}\label{subsec_examples}

In contrast to the Chas-Sullivan product there is not a big range of examples where the string topology coproduct has been computed.
In this subsection we try to list all the examples where a (partial) computation has been carried out.

As already mentioned Naef carries out some partial computations for Lens spaces in \cite{naef:2021}.
Indeed there are homology classes in the free loop spaces of Lens spaces with non-trivial string topology coproduct.

In \cite{goresky:2009} Goresky and Hingston carry out partial computations of the Goresky-Hingston product on manifolds admitting a metric where all geodesics are closed and of the same length.
They use mostly Morse-theoretic methods which in the case of spheres can be used to give a full computation of the Goresky-Hingston product.
Note that the standard examples of manifolds admitting a metric where all geodesics are closed and of the same length are the compact rank one symmetric spaces, i.e. the spheres and projective spaces.
In particular, we want to stress that for all manifolds admitting a metric where all geodesics are closed and of the same length the Goresky-Hingston algebra has non-nilpotent elements.

The knowledge of the Goresky-Hingston product for odd-dimensional spheres was used by Hingston and Wahl in \cite{hingston:2017} to compute the string topology coproduct for odd-dimensional spheres.
Moreover, the second author of this article used Morse-theoretic methods to give a complete computation of the string topology coproduct and the Goresky-Hingston product with rational coefficients for complex and quaternionic projective space, see \cite{stegemeyer:2022}.

For the following theorem we define an integer $\lambda_M$ by 
$$    \lambda_{\mathbb{S}^n} = n-1, \,\,\lambda_{\mathbb{C}P^n} = 1, \,\, \lambda_{\mathbb{H}P^n} = 3 \quad \text{and}\quad \lambda_{\mathbb{O}P^2} = 7.    $$
We sum up these results in the following theorem.
\begin{theorem}\label{theorem_non_trivial_rank_one}
    Let $M$ be a sphere or a projective space.
    Then the string topology coproduct and the Goresky-Hingston product are non-trivial.
    In particular, there is a non-nilpotent cohomology class $\omega\in \mathrm{H}^{\lambda_M}(\Lambda M,M)$.
\end{theorem}
The powers of the cohomology class $\omega$ under the Goresky-Hingston product are related to the iterations of closed geodesics with respect to the symmetric metrics on the compact symmetric spaces of rank one.
We shall revisit this fact in the next subsection.

Based on the results in \cite{goresky:2009}, Hingston partially computes the Goresky-Hingston product for the connected sum $\mathbb{C}P^n\# \mathbb{C}P^n$ in \cite{hingston:2010} and shows that there are non-trivial products.

As far as the authors are aware these are all the examples of spaces where a non-trivial string topology coproduct has been computed.

In contrast to Theorem \ref{theorem_non_trivial_rank_one} the second author has shown the following in \cite{stegemeyer:2021}.
\begin{theorem}[Theorem 5.6 in \cite{stegemeyer:2021}]
    Let $G$ be a simply connected compact Lie group of rank $r\geq 2$.
    Then the string topology coproduct on $G$ is trivial with arbitrary coefficients.
\end{theorem}
This is shown by using the notion of \textit{intersection multiplicity} which we shall introduce in the next Section.
We shall see in Theorem \ref{theorem_torus} that an analogous results holds for tori.

Note that both compact Lie groups as well as spheres and projective spaces are compact symmetric spaces.
Hence, it seems natural to wonder how the string topology coproduct behaves for compact symmetric spaces of higher rank.
The authors of this article study this in \cite{kupper:2022} using the Morse theory of the energy functional of the symmetric metric.
In particular, for each simply-connected compact symmetric space $M$ there are homology classes $W_M^k\in \mathrm{H}_{\bullet}(\Lambda M,M;\mathbb{Z}_2)$, $k\in \mathbb{N}$, which satisfy the following:
\begin{itemize}
    \item For the compact symmetric spaces of rank one these classes have non-trivial coproduct and they are dual to the non-nilpotent elements of the Goresky-Hingston algebra under the duality of homology and cohomology.
    \item If $M$ is a simply-connected compact symmetric space of rank $r\geq 2$ then it holds that $\vee W_M^k = 0$.
\end{itemize}
See \cite[Section 6]{kupper:2022} for a more precise statement and more details.
In particular one sees that the Goresky-Hingston product behaves very differently when one compares symmetric spaces of rank one and symmetric spaces of higher rank.

In fact the authors conjecture that the string topology coproduct is trivial for all irreducible compact symmetric spaces of rank $r\geq 2$ with Euler characteristic $0$. 

Another triviality result for the string topology coproduct was achieved by Naito in \cite{naito:2021}.
He shows that the string topology coproduct with rational coefficients is trivial for a large class of spaces, which in particular includes products of odd-dimensional spheres.
\begin{theorem}[Theorem 1.0.1 in \cite{naito:2021}]
    Let $M$ be a closed simply connected pure manifold.
    If $\mathrm{dim}(\pi_{\mathrm{odd}}(M)\otimes \mathbb{Q}) - \mathrm{dim}(\pi_{\mathrm{even}}(M)\otimes \mathbb{Q})\geq 2$ then the rational string topology coproduct is trivial.
\end{theorem}
For the definition and properties of \textit{pure manifolds}, see \cite[Section 6.1]{naito:2021}.

In Section \ref{sec_results} we shall see more classes of examples with trivial string topology coproduct or at least with trivial based string topology coproduct.
In particular, we will show that for all product manifolds $M = N\times P$ with both $N$ and $P$ having vanishing Euler characteristic the string topology coproduct on $M$ vanishes.

It seems to the authors that there are only very few spaces where it is known that the string topology coproduct is non-trivial while we have large families of examples of manifolds with trivial string topology coproduct.
We therefore ask the following questions.
\begin{question}\hphantom{bs}

    Are there more examples of manifolds apart from the ones we mentioned above where the string topology coproduct is non-trivial? 
    
    Is it possible to characterize manifolds which have a non-trivial string topology coproduct?
\end{question}

\subsection{Geometric Significance}\label{subsec_geometry}

We conclude this section by citing some results which relate the string topology operations, in particular the string topology coproduct and the Goresky-Hingston product, to the theory of closed geodesics in a Riemannian manifold.

Let $(M,g)$ be a closed Riemannian manifold.
Then a \textit{closed geodesic} $c\in \Lambda M$ is a geodesic $c\colon I\to M$ such that
$$     c(0) = c(1) \quad \text{and}\quad \Dot{c}(0) = \Dot{c}(1) .         $$
It is well-known that the closed geodesics in $M$ with respect to the metric $g$ are the critical points of the energy functional $\mathcal{L}_g\colon \Lambda M\to \mathbb{R}$ given by
$$   \mathcal{L}_g(\gamma) = \sqrt{\int_0^1 g_{\gamma(t)}(\Dot{\gamma}(t),\Dot{\gamma}(t)) \,\mathrm{d}t}  . $$
The function $\mathcal{L}_g$ is a smooth function on $\Lambda M\setminus M$ and $c\in \Lambda M$ is a closed geodesic if and only if $(\mathrm{d}\mathcal{L}_g)_c = 0$. Morse theory gives a link between properties of the closed geodesics in $M$ and the homology of the free loop space.
We refer to \cite{klingenberg:78} and \cite{klingenberg:1995} for an introduction to closed geodesics and the relation of this theory to the free loop space.
For a recent survey of the Morse theoretic treatment of closed geodesics we refer to \cite{oancea:2015}.

It is now a natural question if we can relate the string topology operations to properties of closed geodesics and vice versa.
Note that this question is one of the main motivations behind the construction of the Goresky-Hingston product itself in \cite{goresky:2009}.
One Morse-theoretic way of relating homology classes to closed geodesics is the notion of a \textit{critical value} of a homology class.

Before we state the definition, we want to recall some basic Morse-theoretic concepts.
Let $f\colon X\to \mathbb{R}$ be a smooth function on a Hilbert manifold $X$.
Then we denote the \textit{sublevel sets} of $X$ with respect to $f$ by
$$      X^{\leq a} = \{ x\in X\,|\, f(x)\leq a\} \quad \text{and}\quad X^{<a} = \{x\in X\,|\, f(x)<a\}      $$
for $a\in \mathbb{R}$.
Moreover, we recall that in the situation of infinite-dimensional Hilbert manifolds one often asks for the respective functions to satisfy condition (C), see \cite[Section 1.4]{klingenberg:78}.

\begin{definition}[\cite{hingston:2013}, p. 142] \label{cr_val_hom_class}
Let $X$ be a Hilbert manifold and assume that $f:X\to \mathbb{R}$ is a smooth function satisfying condition (C).
\begin{enumerate}
    \item 
Let $b\in\mathbb{R} $ and let $A\in \mathrm{H}_{\bullet}(X,X^{\leq b})$ be a non-trivial homology class where homology is taken with arbitrary coefficients.
Then the \textit{critical value} of $A$ is defined by
$$   \mathrm{cr}(A) = \mathrm{inf}\{  a\in [b,\infty)\,|\, A\in \mathrm{im}(\mathrm{H}_{\bullet}(X^{\leq a},X^{\leq b}) \to \mathrm{H}_{\bullet}(X, X^{\leq b}) )\}     $$
where the morphism is induced by the inclusion $(X^{\leq a}, X^{\leq b}) \hookrightarrow (X, X^{\leq b})$.
\item Let $b\in \mathbb{R}$ and let $\xi\in \mathrm{H}^{\bullet}(X,X^{\leq b})$ be a non-trivial cohomology class where cohomology is taken with arbitrary coefficients.
Then the \textit{critical value} of $\xi$ is defined by
$$    \mathrm{cr}(\xi) = \mathrm{sup}\{ a\in(b,\infty)\,|\,\xi\in \mathrm{im}(       \mathrm{H}^{\bullet}(X,X^{<a})  \to \mathrm{H}^{\bullet}(X,X^{\leq b})    )\}      $$
where the morphism is induced by the inclusion $(X,X^{\leq b}) \hookrightarrow (X,X^{<a})$.
\end{enumerate}
\end{definition}
It turns out that the critical value of a homology class or a cohomology class is indeed a critical value of the function $f$. See also \cite[Section 4]{goresky:2009} for some properties of critical values of homology classes.

In the following we shall now fix a Riemannian metric $g$ on $M$ and we shall consider the critical values with respect to the energy functional $\mathcal{L}_g$.
The fundamental properties of the Chas-Sullivan product and the Goresky-Hingston product with respect to the critical values are the inequalities
\begin{equation}\label{eq_cr_values}
         \mathrm{cr}(A\wedge B) \leq \mathrm{cr}(A) + \mathrm{cr}(B)  \quad \text{and}\quad \mathrm{cr}(\varphi\ostar\psi) \geq \mathrm{cr}(\varphi) + \mathrm{cr}(\psi)    
\end{equation}
for homology classes $A,B\in\mathrm{H}_{\bullet}(\Lambda M)$ and cohomology classes $\varphi,\psi\in\mathrm{H}^{\bullet}(\Lambda M,M)$.
We refer to \cite{goresky:2009} for a proof.

As we have mentioned in the previous subsection, Goresky and Hingston show in \cite{goresky:2009} that the Goresky-Hingston product and the Chas-Sullivan product on compact rank one symmetric spaces can be related to the iteration of closed geodesics.
Here, by iteration we mean the fact that for $m\in \mathbb{N}$ and a closed geodesic $c\in \Lambda M$ the loop $c^m(t) = c(mt)$ is again a closed geodesic.

Further results relating string topology to geometric properties of closed geodesics were achieved by Hingston and Rademacher in \cite{hingston:2013}.
\begin{theorem}[Theorem 1.1 and Corollary 1.3 in \cite{hingston:2013}]
Consider a sphere $\mathbb{S}^n$ of dimension $n\geq 3$. Take a fixed Riemannian metric $g$ and a coefficient field $\mathbb{K}$.
\begin{enumerate}
    \item 
There is a constant $\alpha_{g,\mathbb{K}}  $ such that the difference
$$    \mathrm{deg}(X) - \alpha_{g,\mathbb{K}} \, \mathrm{cr}(X)    $$
is bounded as $X$ runs over all non-trivial homology classes $X\in\mathrm{H}_{\bullet}(\Lambda M;\mathbb{K})$.
\item Let $n\geq 3$ be odd and take $\mathbb{K}= \mathbb{Q}$. Assume that the metric $g$ has sectional curvature satisfying $\tfrac{1}{4}< K \leq 1$.
Then there are at least two closed geodesics with mean frequency equal to $\alpha_{g,\mathbb{Q}}$ or there is a sequence of distinct prime closed geodesics with mean frequencies converging to $\alpha_{g,\mathbb{Q}}$.
\end{enumerate}
\end{theorem}
For the notion of \textit{mean frequency} of a closed geodesic we refer to \cite{hingston:2013}.
We want to stress that the result by Hingston and Rademacher is achieved by the explicit computation of both the Chas-Sullivan ring as well as the Goresky-Hingston ring for spheres.
Moreover they also use the $\mathbb{S}^1$-action of the free loop space and a very precise analysis of the inequalities \eqref{eq_cr_values}.

Finally, we want to mention a result by the second author that uses the non-triviality of the string topology coproduct to detect high multiplicities of closed geodesics.
We say that a closed geodesic $c\in \Lambda M$ is \textit{prime} if it is not the iterate of a closed geodesic $\sigma$.
We say that a closed geodesic $c$ has \textit{multiplicity} $m$ if it is the $m$-th iterate of a prime closed geodesic $\sigma$, i.e. $c = \sigma^m$.
Note that there are iterated versions of the coproduct, see \cite[Section 3.3]{hingston:2017}, which are maps
$$    \vee^k\colon \mathrm{H}_i(\Lambda ,M) \to \mathrm{H}_{i+k-kn}(\Lambda^{k+1}, \Lambda^k \times M\cup \Lambda^{k-1}\times M\times \Lambda \cup \ldots \cup M\times \Lambda^k)   .   $$
If we take $k = 1$ then this gives the ordinary string topology coproduct and if we take coefficients in a field then we can understand $\vee^k$ as the composition
$$        (\underbrace{\mathrm{id}\otimes\ldots\otimes \mathrm{id}}_{k-1\,\,\text{times}}\otimes {\vee})\circ (\underbrace{\mathrm{id}\otimes\ldots\otimes \mathrm{id}}_{k-2\,\,\text{times}}\otimes {\vee})\circ \ldots \circ {\vee} : \mathrm{H}_{\bullet}(\Lambda ,M)\to (\mathrm{H}_{\bullet}(\Lambda ,M))^{\otimes k+1}  .  $$
On the homology of the $k$-fold product of $\Lambda M$ with itself, $k\in\mathbb{N}$, one can also define critical values of homology classes with respect to the functional
$$     \mathcal{L}_{g,k}\colon \Lambda M^k \to \mathbb{R} , \quad  \mathcal{L}_{g,k}(\gamma_1,\ldots,\gamma_k) = \sum_{i=1}^k \mathcal{L}_g(\gamma_i) ,   $$
see \cite{stegemeyer2023global}.
The second author of this article then shows the following.
\begin{theorem}[Theorem 2.3.8 in \cite{stegemeyer2023global}]
    Let $R$ be a commutative unital ring and take homology with coefficients in $R$.
    Let $(M,g)$ be a closed $R$-oriented Riemannian manifold of dimension $n$ and assume that the energy functional $\mathcal{L}_g$ is a Morse-Bott function.
Let $X\in\mathrm{H}_{\bullet}(\Lambda M,M)$ be a homology class with critical value $\mathrm{cr}(X)= a>0$ and let $\Sigma^a\subseteq \Lambda M$ be the critical set at level $a$.
If the critical value of the $k$-fold iterated coproduct $\vee^k X$ satisfies $\mathrm{cr}(\vee^k X) = a = \mathrm{cr}(X)$, then there is a closed geodesic $c \in \Sigma^a$ with multiplicity greater or equal to $k+1$.
\end{theorem}

We remark that on closed manifolds being Morse-Bott is a generic property of the energy functional of Riemannian metrics, see \cite{abraham:1970}.

\section{Intersection Multiplicity and Triviality of the Coproduct}\label{sec_int_mult}

In this section we turn to the concept of intersection multiplicity.
We begin by introducing the notion of intersection multiplicity of homology classes in loop spaces and will then see the connection to the string topology coproduct.

In \cite[Section 5]{hingston:2017} Hingston and Wahl define the \textit{basepoint intersection multiplicity} of a homology class in the free loop space.
We shall also define an analogous notion for the based loop space.

\begin{definition} Let $M$ be a closed manifold with basepoint $p_0\in M$ and let $R$ be a commutative unital ring.
\begin{enumerate}
    \item 
        Let $[X]\in \mathrm{H}_{*}(\Lambda M,M)$ be a homology class represented by a cycle $X\in \mathrm{C}_{*}(\Lambda M,M)$.
         Assume that the relative cycle $X$ is represented by a cycle $x\in \mathrm{C}_{*}(\Lambda M)$.
        The \textit{intersection multiplicity} $\mathrm{int}_{\Lambda}([X])$ of the homology class $[X]$ is the integer
     $$    \mathrm{int}_{\Lambda}([X]) = \inf_{A \sim x} \Big(\sup\big[ \#(\gamma^{-1}(\{p_0\})) \,|\, \gamma\in \mathrm{im}(A),\,\,\, \mathcal{L}_g(\gamma)> 0 \big] \Big)  -1  . $$
    The infimum is taken over all cycles $A \in \mathrm{C}_{*}(\Lambda M)$ homologous to $x$ and $\mathcal{L}_g$ is the energy functional for an arbitrary Riemannian metric $g$.
    \item
     Let $[X]\in \mathrm{H}_{*}(\Omega M,\{p_0\})$ be a homology class with representing cycle $X\in \mathrm{C}_{*}(\Omega M,\{p_0\})$.
         Assume that the relative cycle $X$ is represented by a cycle $x\in \mathrm{C}_{*}(\Omega M)$.
        The \textit{intersection multiplicity} $\mathrm{int}_{\Omega}([X])$ of the homology class $[X]$ is the integer
     $$    \mathrm{int}_{\Omega}([X]) = \inf_{A \sim x} \Big(\sup\big[ \#(\gamma^{-1}(\{p_0\})) \,|\, \gamma\in \mathrm{im}(A),\,\,\, \mathcal{L}_g(\gamma)> 0 \big] \Big)  -1  . $$
    The infimum is taken over all cycles $A \in \mathrm{C}_{*}(\Omega M)$ homologous to $x$ and $\mathcal{L}_g$ is the energy functional for an arbitrary Riemannian metric $g$.
\end{enumerate}
\end{definition}

As Hingston and Wahl remark in \cite[Section 5]{hingston:2017} the intersection multiplicity of a homology class is always finite and independent of the choice of the metric $g$.
Hingston and Wahl study properties of the string topology coproduct that can be deduced from knowing the intersection multiplicity of a homology class.
For details on this and for a proof of the next theorem we refer to \cite{hingston:2017}.
\begin{theorem}[\cite{hingston:2017}, Theorem 3.10 and Theorem 5.2]\label{int_mult_copro_vanish}
    Let $M$ be a closed manifold. Take homology with coefficients in a commutative unital ring $R$ and assume that $M$ is $R$-oriented.
    \begin{enumerate}
        \item If $[X]\in\mathrm{H}_{\bullet}(\Lambda M,M)$ is a homology class in the free loop space with intersection multiplicity $\mathrm{int}_{\Lambda}([X])\leq k$, then $\vee^k[X] = 0$.
        \item If $M$ is a sphere or a projective space then the converse also holds, i.e. if $[X]$ satisfies $\vee^k[X] = 0$ then $\mathrm{int}_{\Lambda}([X])\leq k$.
    \end{enumerate}
\end{theorem}

Since the proof of \cite[Theorem 3.10]{hingston:2017} which includes Lemma 3.5 and Proposition 3.6 of \cite{hingston:2017} only uses maps that keep the basepoints of loops fixed, the proof can be transferred directly to the case of the based loop space.
One obtains the following result.
\begin{proposition}[\cite{hingston:2017}, Theorem 3.10] \label{prop_intersection_trivial}
Let $M$ be a closed manifold with basepoint $p_0\in M$ and let $R$ be a commutative unital ring.
Assume that $M$ is $R$-oriented.
 Suppose that $[X]\in \mathrm{H}_{\bullet}(\Omega_{p_0}M,\{p_0\})$ has intersection multiplicity $\mathrm{int}_{\Omega}([X]) = 1$, then the based coproduct $\vee_{\Omega}[X]$ vanishes.
\end{proposition}

In the following theorem we sum up some known results about the intersection multiplicity of particular classes of spaces.

\begin{theorem}
    Let $M$ be a closed manifold.
    \begin{enumerate}
        \item If $M = G$ is a compact simply connected Lie group of rank $r\geq 2$ and if we take homology with coefficients in a commutative unital ring $R$, then every non-trivial homology class $[X]\in\mathrm{H}_{\bullet}(\Lambda G,G)$ has intersection multiplicity $\mathrm{int}_{\Lambda}([X]) = 1$. Consequently, the string topology coproduct on $\mathrm{H}_{\bullet}(\Lambda G,G)$ vanishes.
            \item If $M$ is a compact simply connected symmetric space of rank $r\geq 2$ and if we consider homology with coefficients $R = \mathbb{Z}_2$ then every non-trivial homology class $[X]\in\mathrm{H}_{\bullet}(\Omega M,\{p_0\})$ has intersection multiplicity $\mathrm{int}_{\Omega}([X]) = 1$.
             Consequently, the based string topology coproduct on $\mathrm{H}_{\bullet}(\Omega M,\{p_0\})$ vanishes.
    \end{enumerate}
\end{theorem}
For proofs we refer to \cite{kupper:2022} and \cite{stegemeyer:2021}.
We now want to study some elementary properties of intersection multiplicity.
First, let us note the following obvious properties.
\begin{proposition}\label{prop_inclusion}
    Let $M$ be a closed manifold. Take homology with coefficients in a commutative unital ring $R$ and assume that $M$ is $R$-oriented.
    Let $j\colon \Omega M\hookrightarrow \Lambda M$ be the inclusion.
    \begin{enumerate}
        \item If $[X]\in\mathrm{H}_{\bullet}(\Omega M,\{p_0\})$ satisfies $\mathrm{int}_{\Omega}([X]) \leq k$ then $\mathrm{int}_{\Lambda}(j_* [X]) \leq k$.
        \item If every class $[X]\in\mathrm{H}_{\bullet}(\Omega M,\{p_0\})$ satisfies $\mathrm{int}_{\Omega}([X]) = 1$ then the string topology coproduct vanishes on the image of $j_*\colon \mathrm{H}_{\bullet}(\Omega M,\{p_0\})\to \mathrm{H}_{\bullet}(\Lambda M,M)$.
    \end{enumerate}
\end{proposition}
\begin{proposition}\label{prop_products_one_factor_trivial_int_mult}
Let $M = N\times P$ be a product of two closed manifolds. Take homology with coefficients in a commutative unital ring $R$.
    \begin{enumerate}
        \item Assume that $\mathrm{H}_{\bullet}(\Lambda N)$ or $\mathrm{H}_{\bullet}(\Lambda P)$ is a free $R$-module. If every homology class $[X]\in\mathrm{H}_{\bullet}(\Lambda N)$ can be represented by a cycle $x\in \mathrm{C}_{\bullet}(\Lambda N)$ such that each loop $\gamma\in \mathrm{im}(x)$ has no non-trivial self-intersections then each class $[Z]\in\mathrm{H}_{\bullet}(\Lambda M,M)$ has intersection multiplicity $\mathrm{int}_{\Lambda}([Z])= 1$.        
        \item In the situation as above assume that both $N$ and $P$ are $R$-oriented. Then the string topology coproduct of $M$ vanishes.
    \end{enumerate}
\end{proposition}
\begin{proof}
    Note that we have $\Lambda M = \Lambda (N\times P) \cong \Lambda N\times \Lambda P$.
    By the assumption that one of $\mathrm{H}_{\bullet}(\Lambda N)$ or $\mathrm{H}_{\bullet}(\Lambda P)$ is free, we have a Künneth isomorphism $$\mathrm{H}_{\bullet}(\Lambda N)\otimes \mathrm{H}_{\bullet}( \Lambda P) \cong \mathrm{H}_{\bullet}(\Lambda M). $$
    Hence, the classes of the form $[Z] = [X]\times [Y] \in\mathrm{H}_{\bullet}(\Lambda M)$ with $[X]\in\mathrm{H}_{\bullet}(\Lambda N)$ and $[Y]\in\mathrm{H}_{\bullet}(\Lambda P)$ generate the homology of $\Lambda M$.
    Choose a representing cycle $x\in \mathrm{C}_{\bullet}(\Lambda N)$ of $[X]$ such that every loop in its image has no non-trivial self-intersections.
    Let $y\in\mathrm{C}_{\bullet}(\Lambda P)$ be an arbitrary representing cycles of $[Y]$, then $z=x\times y$ represents $[Z]$.
    Let $\gamma = (\sigma,\eta)\in \mathrm{im}(z)$ be a loop in the image of $z$.
    Note that $\sigma\in\mathrm{im}(x)$, thus by the choice of $x$ the loop $\sigma$ has no non-trivial self-intersections and consequently, $\gamma$ has no non-trivial self-intersections either.
    Since this holds for all loops in the image of $z$ we found a cycle representing $[Z]$ such that all loops in its image have no non-trivial self-intersections.
    If now $N$ and $P$ are $R$-oriented then by Theorem \ref{int_mult_copro_vanish} the coproduct on $N\times P$ vanishes.
\end{proof}
With the considerations of \cite{stegemeyer:2021} we obtain the following Corollary.
\begin{corollary}
    Let $M = G\times P$ be a product of a simply connected Lie group $G$ of rank $r\geq 2$ with a closed manifold $P$. Let $R$ be a commutative unital ring and take homology with $R$-coefficients.
    Assume that $P$ is $R$-oriented and that $\mathrm{H}_{\bullet}(\Lambda G)$ or $\mathrm{H}_{\bullet}(\Lambda P)$ is a free $R$-module.
    Then the string topology coproduct  vanishes on the product $M = G\times P$.
\end{corollary}
\begin{proof}
    There is a splitting 
    $$      \mathrm{H}_{\bullet}(\Lambda G)\cong \mathrm{H}_{\bullet}(\Lambda G,G)\oplus \mathrm{H}_{\bullet}(G) .   $$
    From \cite{stegemeyer:2021} we know that all classes in $\mathrm{H}_{\bullet}(\Lambda G,G)$ can be represented by cycles in $\mathrm{C}_{\bullet}(\Lambda G)$ such that each loop in the image of the respective cycle has no non-trivial self-intersection.
    Moreover, since $G$ is a Lie group it has Euler characteristic $0$.
    Therefore we can represent the homology classes in $\mathrm{H}_{\bullet}(G)\hookrightarrow \mathrm{H}_{\bullet}(\Lambda G)$ coming from the constant loops by a loop with no non-trivial self-intersection.
    Hence the assumptions of Proposition \ref{prop_products_one_factor_trivial_int_mult} are satisfied and the claim follows.
\end{proof}
Now, assume that we have an action of a topological group $\varphi\colon G\times M\to M$ on the manifold $M$.
This action induces a continuous action of $G$ on the free loop space $\Phi\colon G\times \Lambda M\to \Lambda M$ given by
$$    \Phi(g,\gamma)(t) =   \varphi(g,\gamma(t)) .   $$
Since $G$ acts by homeomorphisms it is clear that if $\gamma\in \Lambda M$ has precisely $k$ non-trivial self-intersections, then $\Phi(g,\gamma)$ also has precisely $k$ non-trivial self-intersections.
The group action $\Phi$ induces a pairing $\Theta\colon \mathrm{H}_i(G)\otimes \mathrm{H}_j(\Lambda M,M)\to \mathrm{H}_{i+j}(\Lambda M,M)$ given by
$$    \Theta(X\otimes Y) = \Phi_* (X\times Y)\quad \text{for}\,\,\, X\in \mathrm{H}_i(G), \,\,\, Y\in\mathrm{H}_j(\Lambda M,M).   $$
We note the following properties for intersection multiplicity.
\begin{proposition}\label{prop_group_action}
    Let $M$ be a closed manifold with a group action $\varphi\colon G\times M\to M$.
    Take homology with coefficients in a commutative unital ring $R$ and assume that $M$ is $R$-oriented.
    Let $\Phi\colon G\times \Lambda M\to \Lambda M$ be the induced action on the free loop space.
    \begin{enumerate}
        \item If $[X]\in\mathrm{H}_{\bullet}(\Lambda M,M)$ is a relative homology class with intersection multiplicity $\mathrm{int}_{\Lambda}([X])\leq k$ and $[A]\in\mathrm{H}_{\bullet}(G)$, then $\Theta([A]\otimes [X])$ satisfies $\mathrm{int}_{\Lambda}(\Theta([A]\otimes [X]) \leq k$.
        \item If every homology class $[X]$ in the image of $j_*\colon \mathrm{H}_{\bullet}(\Omega M,\{p_0\})\to \mathrm{H}_{\bullet}(\Lambda M,M)$ has intersection multiplicity $1$ then the string topology coproduct vanishes on the subspace
        $$    U  = \Theta \big( \mathrm{H}_{\bullet}(G)\otimes \mathrm{im}\left(j_*\colon \mathrm{H}_{\bullet}(\Omega M,\{p_0\})\to \mathrm{H}_{\bullet}(\Lambda M,M)\right)\big) .   $$
    \end{enumerate}
\end{proposition}
For the second statement we also refer to \cite[Proposition 6.6]{kupper:2022}.
We shall now use the second statement of the above Proposition to prove that the string topology coproduct vanishes for all tori of rank $r\geq 2$.
\begin{theorem}\label{theorem_torus}
    Consider the $n$-dimensional torus $T^n$ with $n\geq 2$ and take homology with coefficients in a commutative unital ring $R$.
    \begin{enumerate}
        \item Every class $[X]\in\mathrm{H}_{\bullet}(\Omega T^n,\{p_0\})$ satisfies $\mathrm{int}_{\Omega}([X]) = 1$.
        \item Every class $[X]\in\mathrm{H}_{\bullet}(\Lambda T^n,T^n)$ satisfies $\mathrm{int}_{\Lambda}([X]) = 1$.
        \item Both the based string topology coproduct and the string topology coproduct are trivial on $T^n$.
    \end{enumerate}
\end{theorem}
\begin{proof}
    We start by noting that the based loop space of $\mathbb{S}^1$ is homotopy equivalent to the discrete set $\mathbb{Z}$.
    This is because we have $\pi_0(\Omega \mathbb{S}^1)\cong \mathbb{Z}$ and $\pi_i((\Omega \mathbb{S}^1 )_j)= \{0\}$ for each path component $(\Omega \mathbb{S}^1)_j$ and each $i\geq 1$.
    Therefore we have 
    $$     \mathrm{H}_i(\Omega \mathbb{S}^1,\{p_0\}) \cong     \begin{cases} \bigoplus_{k\in\mathbb{Z}\setminus\{0\}} R & i = 0\\ \{0\} & \text{else} .
           
    \end{cases}  $$
    Moreover, the generators are given by the standard $k$-fold covering of the circle, i.e. the non-trivial closed geodesics.
    Since $T^n = (\mathbb{S}^1)^n$ we see that
    $$     \mathrm{H}_i(\Omega T^n,\{p_0\}) \cong     \begin{cases} \bigoplus_{(k_1,\ldots,k_n)\in\mathbb{Z}^n\setminus\{0\}} R & i = 0\\ \{0\} & \text{else} .     
    \end{cases}  $$
    The generators of $\mathrm{H}_0(\Omega T^n,\{p_0\})$ are given by the homology classes generated by the closed geodesics in $T^n$ with respect to the standard flat metric starting at the basepoint and considered as a point in $\Omega T^n$.
    Let $\gamma\in \Omega T^n$ be such a closed geodesic and consider the induced homology class $[\gamma]\in\mathrm{H}_0(\Omega T^n,\{p_0\})$.
    If $\gamma$ is prime then it has no non-trivial intersections of the basepoint and it is clear that $\mathrm{int}_{\Omega}([\gamma]) = 1$.
    If $\gamma$ does intersect the basepoint then we choose a homotopic path $\sigma\in \Omega M$ with no non-trivial self-intersections.
    This is possible since the dimension of the torus is greater or equal to $2$.
    Since $\sigma\simeq \gamma$ they represent the same homology class, i.e. $[\sigma] = [\gamma]$.
    Since the loop $\sigma$ does not have any non-trivial self-intersections we obtain $\mathrm{int}_{\Omega}([\gamma]) = 1$.
    Since all non-trivial homology classes in $\mathrm{H}_{\bullet}(\Omega T^n,\{p_0\})$ are of this form this shows the first part.
    For the second part note that $T^n$ is a Lie group and therefore the map
    $$    F\colon T^n\times \Omega T^n\to \Lambda T^n ,\quad F(x,\gamma)(t) = x + \gamma(t)       $$
    is a homeomorphism.
    Moreover, both $\mathrm{H}_{\bullet}(T^n)$ and $\mathrm{H}_{\bullet}(\Omega T^n,\{p_0\})$ are free and thus the pairing
    $$    \mathrm{H}_{\bullet}(T^n) \otimes \mathrm{H}_{\bullet}(\Omega T^n,\{p_0\}) \xrightarrow[]{\times }\mathrm{H}_{\bullet}(T^n\times \Omega T^n, T^n\times \{p_0\}) \xrightarrow[]{F_*} \mathrm{H}_{\bullet}(\Lambda T^n,T^n)  $$
    is an isomorphism.
    This pairing is precisely the pairing induced by the group action of $T^n$ on itself, see the discussion before Proposition \ref{prop_group_action}.
    Therefore, by Propositions \ref{prop_inclusion} and \ref{prop_group_action} we see that every class $[X]\in\mathrm{H}_{\bullet}(\Lambda T^n,T^n)$ has intersection multiplcitity $\mathrm{int}_{\Lambda}([X]) = 1$.
    This shows the second part and the third part then follows immediately.
\end{proof}

\section{Results on intersection multiplicity}\label{sec_results}

In this section we prove several results on intersection multiplicity both for the free loop space and for the based loop space.
We also discuss the implications for the string topology coproduct.

Before we start we want to remark that in order to show that all homology classes in the relative homology $\mathrm{H}_{\bullet}(\Lambda M,M)$ for $M$ a closed manifold have intersection multiplicity $1$ it is sufficient to show that all homology classes in the absolute homology $\mathrm{H}_{\bullet}(\Lambda M)$ can be represented by a cycle such that all loops in its image have no non-trivial self-intersections.
This is because we have that the inclusion $i\colon M\to \Lambda M$ is a section of the free loop fibration $\mathrm{ev}\colon \Lambda M\to M$ and therefore
$$    \mathrm{H}_i(\Lambda M) \cong \mathrm{H}_i(M)\oplus \mathrm{H}_i(\Lambda M,M) \quad \text{for all}\,\,\,i\in\mathbb{N}_0 .    $$
The analogous property holds for the based loop space.

\subsection{Product manifolds}
Let $M$ and $N$ be compact manifolds with vanishing Euler characteristic.
Then there exist nowhere vanishing vector fields $X$ on $M$ and $Y$ on $N$.
Let $\Phi\colon M\times \mathbb{R}\to M$ and $\Psi\colon N\times \mathbb{R}\to N$ be the corresponding flows.
By potentially rescaling $X$ and $Y$ we can achieve that $\Phi$ and $\Psi$ satisfy
$$      \Phi(p,s) \neq p\quad \text{for all}\,\,\,p\in M, s\in (0,1]       $$
and
$$      \Psi(q,s) \neq q\quad \text{for all}\,\,\, q\in N, s\in (0,1]  .     $$
We shall use these flows to deform cycles on the free loop space of $M\times N$ such that no loop in the image of these cycles has a basepoint intersection.
\begin{theorem}\label{theorem_int_mult_products}
    Let $R$ be a commutative unital ring and take homology with coefficients in $R$.
    Let $M$ and $N$ be closed manifolds both of Euler characteristic $0$.
    Then every class $Z\in\mathrm{H}_{\bullet}(\Lambda (M\times N),M\times N)$ has intersection multiplicity $\mathrm{int}_{\Lambda}(Z) = 1$.
\end{theorem}
\begin{proof}
Let $\Phi\colon M\times \mathbb{R}\to M$ and $\Psi\colon N\times \mathbb{R}\to N$ 
 be the flows as above.
 We define maps
$    f\colon \Lambda M\to \Lambda M $ and $g\colon \Lambda N\to \Lambda N$ by setting
$$   f(\sigma)(t) = \begin{cases}
    \Phi(\sigma(0),3t), & 0\leq t\leq \tfrac{1}{3} \\
    \Phi(\sigma(0), 2-3t), & \tfrac{1}{3}\leq t\leq \tfrac{2}{3} \\
    \sigma(3t-2), & \tfrac{2}{3}\leq t\leq 1 
\end{cases}    $$
for $\sigma\in \Lambda M$ and
$$   g(\eta)(t) = \begin{cases}
    \eta(3t) , & 0\leq t\leq \tfrac{1}{3} \\
    \Psi(\eta(1),3t-1), & \tfrac{1}{3} \leq t\leq \frac{2}{3} \\
    \Psi(\eta(1),3-3t), & \tfrac{2}{3} \leq t \leq 1
\end{cases}     $$
for $\eta\in \Lambda N$.
One checks that $f\simeq \mathrm{id}_{\Lambda M}$ and $g\simeq \mathrm{id}_{\Lambda N}$ and hence $(f\times g)\simeq \mathrm{id}_{\Lambda (M\times N)}$.
Now, let $Z\in\mathrm{H}_{\bullet}(\Lambda (M\times N))$ be a homology class.
Choose a representing cycle $z'\in\mathrm{C}_{\bullet}(\Lambda (M\times N))$ of $Z$.
We define the chain $z= (f\times g)_* (z')$.
Then $z$ and $z'$ are homologous and thus $z$ represents the class $Z$ as well.
Let $\gamma\in \mathrm{im}(z)$, i.e.
$$   \gamma = (f(\sigma),g(\eta))   $$
for some $\sigma\in \Lambda M$ and $\eta\in \Lambda N$.
Assume that there is a $t_*\in (0,1)$ with $\gamma(0) = \gamma(t_*)$.
This implies that both
$$     f(\sigma)(t_* ) = f(\sigma)(0) = \sigma(0) \quad \text{and} \quad g(\eta)(t_*) = g(\eta)(0)  = \eta(0)  .      $$
By construction of $f$ and $g$ the first condition impies that $t_*\in [\tfrac{2}{3},1)$ while the second condition implies $t_*\in (0,\tfrac{1}{3}]$.
This is a contradiction so we see that $Z$ has intersection multiplicity $1$.
Thus we have shown that all non-trivial homology classes have intersection multiplicity $1$.
\end{proof}
Again, by Theorem \ref{int_mult_copro_vanish} we can infer the following.
\begin{corollary}
    Let $R$ be a commutative unital ring and take homology with coefficients in $R$.
    Let $M$ and $N$ be compact $R$-oriented manifolds both with vanishing Euler characteristic.
    Then the string topology coproduct on $\mathrm{H}_{\bullet}(\Lambda(M\times N),M\times N)$ vanishes.
\end{corollary}
\begin{corollary}
    Let $M = \mathbb{S}^{2m+1}\times \mathbb{S}^{2n+1}$ be a product of odd-dimensional spheres.
    Then the coproduct on $\mathrm{H}_{\bullet}(\Lambda M,M)$ vanishes.
\end{corollary}
\begin{example}
    Theorem \ref{theorem_int_mult_products} gives also an alternative proof that the string topology coproduct vanishes on the torus $T^n$ for $n\geq 2$ since $T^n = \mathbb{S}^1\times T^{n-1}$ and both $\mathbb{S}^1$ as well as $T^{n-1}$ have vanishing Euler characteristic.
\end{example}

We can show a similar result for the based loop space.
\begin{theorem}\label{theorem_int_mult_based}
    Let $R$ be a commutative unital ring and take homology with coefficients in $R$.
     Let $M$ and $N$ be closed manifolds of positive dimension.
    \begin{enumerate}
        \item Every class $Z\in\mathrm{H}_{\bullet}(\Omega (M\times N),\{p_0\})$ has intersection multiplicity $\mathrm{int}_{\Omega}(Z) = 1$.
        \item Assume that $M$ and $N$ are $R$-oriented. Then the based string topology coproduct vanishes on $\mathrm{H}_{\bullet}(\Omega (M\times N),\{p_0\})$.
    \end{enumerate}
\end{theorem}
\begin{proof}
Denote the basepoints of $M$ and $N$ by $p_0\in M$ and $q_0\in N$.
Choose paths $\omega\colon I\to M$ and $\sigma\colon I\to N$ such that 
$$   \omega(0) = p_0, \quad \text{and} \quad \omega(s) \neq p_0 \,\,\,\text{for all}\,\,\, s > 0    $$
and similarly
$$  \sigma(0) = q_0   \quad \text{and} \quad \sigma(s) \neq q_0 \,\,\,\text{for all}\,\,\, s > 0  . $$
Define maps $f\colon \Omega M\to \Omega M$ and $g\colon \Omega N\to \Omega N$ by setting
$$    f(\gamma)(t) = \begin{cases}
        \omega(3t) , & 0\leq t\leq \tfrac{1}{3}
        \\ \omega(2- 3t), & \tfrac{1}{3} \leq 0 \leq \tfrac{2}{3} \\
        \gamma(3t -2), & \tfrac{2}{3} \leq t \leq 1 
\end{cases}
$$
where $\gamma\in \Omega_{p_0 M}$ and
$$
    g(\eta)(t)    = \begin{cases}
        \eta(3t) , & 0\leq t\leq \tfrac{1}{3} 
        \\ \sigma(3t-1), & \tfrac{1}{3} \leq t \leq \tfrac{2}{3} \\
        \sigma(3-3t), & \tfrac{2}{3} \leq t \leq 1
    \end{cases}
$$
for $\eta\in \Omega_{q_0}N$.
Then we have $f\simeq \mathrm{id}_{\Omega M}$ and $g\simeq \mathrm{id}_{\Omega N}$.
As in the proof of Theorem \ref{theorem_int_mult_products} one can now argue that all classes of the form $Z\in\mathrm{H}_{\bullet}(\Omega (M\times N),\{p_0\})$ have intersection multiplicty $1$ since all loops in the image of $f\times g$ have no non-trivial self-intersections.
The second assertion now follows using Proposition \ref{prop_intersection_trivial}.
\end{proof}

\begin{corollary}
    Let $M= \mathbb{S}^m\times \mathbb{S}^n$ be a product of spheres.
    Then the based string topology coproduct of $M$ vanishes.
\end{corollary}

\subsection{Fiber bundle with a section}

We now turn to the situation of a fiber bundle $p\colon E\to B$ with a section.
It turns out that in this case the based loop space $\Omega E$ can be understood particularly well.

Let $p\colon E\to B$ be a smooth fiber bundle with $B$ a compact and connected smooth manifold. We choose a base point $b_0\in B$ and we assume that the fiber $F = p^{-1}(b_0)$ is also a compact and connected smooth manifold. Hence also $E$ is compact and connected. 
Note that it follows that $p\colon E\to B$ is a fibration. 
Let $i\colon F\to E$ denote the fiber inclusion. We furthermore assume that there is a section $s\colon B\to E$ and we choose $s(b_0)$ as a base point in $E$ and $F$. We consider the map
\begin{equation}\label{eq_definition_g}
       g\colon \Omega F\times \Omega B\to \Omega E,\quad (\gamma,\eta)\mapsto \mathrm{concat}(i\circ \gamma,s\circ \eta) \, ,
\end{equation}
where the loops are based at $s(b_0)$ and $b_0$, respectively. 
We now choose $c_{b_0}$ and $c_{s(b_0)}$, the constant loops, as base points in the loop spaces $\Omega B$ and $\Omega E$ as well as $\Omega F$, respectively. It follows that $\Omega s\colon \Omega B\to \Omega E,\,\Omega s(\gamma) =  s\circ \gamma,$ maps $c_{b_0}$ to $c_{s(b_0)}$. It turns out that the homomorphism which is given by the composition
\begin{align*}
        \pi_i(\Omega F\times \Omega B,(c_{s(b_0)},c_{b_0})) \xrightarrow[\cong]{(\mathrm{pr}_1)_* \oplus (\mathrm{pr}_2)_*} \pi_i(\Omega F,c_{s(b_0)})&\oplus \pi_i(\Omega B,c_{b_0})\\
        &\xrightarrow{(\Omega i)_* + (\Omega s)_*} \pi_i(\Omega E,c_{s(b_0)}) 
\end{align*}
equals the map $g_*$ induced by $g$. This holds since the group operation in $\pi_n(\Omega E)$ is defined via the H-space structure of $\Omega E$: Let $x\colon \mathbb{S}^i\to \Omega F$ and $y\colon \mathbb{S}^i\to \Omega B$ be basepoint preserving maps. Then the composition above is given by
$$
[(x,y)]\mapsto [i\circ x]+[s\circ y]
=[p\mapsto \mathrm{concat}(i\circ x(p),s\circ y(p))]
$$
which is exactly $g_*([(x,y)])$. Because $\Omega s$ is a section of the fibration $\Omega p\colon \Omega E\to\Omega B$ we have that the long exact homotopy sequence of $\Omega p$ splits into short exact sequences and exhibits the maps
$$
\pi_i(\Omega F,c_{s(b_0)})\oplus \pi_i(\Omega B,c_{b_0})\xrightarrow{(\Omega i)_* + (\Omega s)_*} \pi_i(\Omega E,c_{s(b_0)})
$$
as isomorphisms for all $i$. 
For $i=0$ replace $+$ by and $\cdot$ and isomorphism by bijection. Hence 
$$
g_*\colon \pi_i\left(\Omega F\times\Omega B,(c_{s(b_0)},c_{b_0})\right)\to\pi_i\left(\Omega E,c_{s(b_0)}\right)
$$
is an isomophism for all $i\geq 1$ and a bijection for $i = 0$. In case that $B,E$ and $F$ are simply-connected \cite[Theorem VII.7.2 (e)]{bredon:2013} now immediately yields that
$$
g_*\colon \pi_n\left(\Omega F\times\Omega B,(\alpha,\beta)\right)\to\pi_n\left(\Omega E,g(\alpha,\beta)\right)
$$
is an isomorphism for arbitrary base points $\alpha\in \Omega F$ and $\beta\in \Omega B$. Consequently, $g$ is a weak homotopy equivalence. Since the loop space of a manifold is a CW complex, see \cite{milnor1959spaces}, the map $g$ is a homotopy equivalence by the Whitehead theorem.

With some more effort one can also show that $g$ is also a homotpy equivalence if the assumption of simply-connectedness is removed. In fact we have
\begin{lemma} \label{lemma_principal_fibration}
    Let $p\colon E\to B$ be a smooth fiber bundle with fiber $F$  and assume that  $B,E,F$ are smooth closed connected manifolds. Also assume that there is a smooth section $s\colon B\to E$. Then the map $g\colon \Omega F\times \Omega B\to \Omega E$ as in equation \eqref{eq_definition_g} is a fiber homotopy equivalence between the trivial fibration $\Omega F\times \Omega B\to \Omega B$ and the fibration $\Omega p\colon\Omega E\to \Omega B$.
\end{lemma}

We remark that $\Omega p\colon\Omega E\to \Omega B$ is a principal fibration and hence trivial up to homotopy if it has a section. 

We use this homotopy equivalence $g$ to prove the following result.

\begin{theorem}\label{thm_fiber_bundle_with_section}
    Let $R$ be a commutative unital ring and take homology with coefficients in $R$.
    Let $p\colon E\to B$ be a smooth fiber bundle with a section, where $B$ and $E$ are closed manifolds. Assume that both $B$ and the fiber $F$ are positive dimensional.
    Then all homology classes in $\mathrm{H}_{\bullet}(\Omega E,\{p_0\})$ have intersection multiplicity $1$.
 \end{theorem}
 
 \begin{proof}
     Let $g\colon \Omega F\times \Omega B\to \Omega E$ be the homotopy equivalence defined in equation \eqref{eq_definition_g}.
      Choose an open neighborhood $U$ of the basepoint $p_0\in E$ such that $U$ is diffeomorphic to a product neighborhood $V= V_1\times V_2$ of $0$ in $\mathbb{R}^k\times \mathbb{R}^m$ and such that $F\cap U$ is taken to $V_1\times \{0\}$ and $s(B)\cap U$ is taken to $\{0\}\times V_2$.
    We can assume that the open cube $(-4,4)^{k+m}\subseteq V$.
    Let $X$ be the vector field on $V$ given by
    $$    X = \sum_{i=1}^n \partial_i     $$
    where $n = k+m$.
    Choose a smooth function $\rho$ which takes the constant value $1$ on $(-2,2)^n$ and vanishes outside $(-3,3)^n$.
    Then let $Y = \rho X$ and extend this to a global vector field on $E$ by continuing it trivially outside $U$.
    Now, let $\Phi\colon E\times \mathbb{R}\to E$ be the flow of this vector field.
    We note that 
    $$   \Phi(p_0,s) = (s,s,\ldots ,s)   \quad \text{for}\,\,\, s\in [-1,1]    $$
    under the identification $U\cong V_1\times V_2\subseteq \mathbb{R}^n$.
    In particular, 
    \begin{equation}\label{eq_property_cross_of_fiber_and_base}
           \Phi(p_0,s)\in V_1\times \{0\} \cup \{0\}\times V_2 \,\,\, \text{for} \,\, s\in [-1,1] \quad \iff \quad s = 0 . 
    \end{equation}
    Define a map $\theta\colon \Omega E\to \Omega E$ by setting
    $$   \theta(\gamma)(s) = \begin{cases}
        \Phi(\gamma(s),2s), & 0\leq s \leq \tfrac{1}{2} \\
        \Phi(\gamma(s), 2(1-s))  & \tfrac{1}{2}\leq s\leq 1 .
    \end{cases}    $$
    One checks that $\theta$ is homotopic to the identity on $\Omega E$.
    Consequently, the composition $\varphi = \theta \circ g\colon \Omega F\times \Omega B\to \Omega E$ is homotopic to $g$.
    Let $X\in \mathrm{H}_{\bullet}(\Omega F\times \Omega B)$ be a non-trivial homology class and let $x\in\mathrm{C}_{\bullet}(\Omega F\times \Omega B)$ be a representing cycle.
    Then the cycle $y = \varphi_{\#} x \in \mathrm{C}_{\bullet}(\Omega E)$ represents the class $g_* X$ and all classes in $\mathrm{H}_{\bullet}(\Omega E)$ have representing cycles of this form since $g$ is a homotopy equivalence.
    It is then clear that $\mathrm{im}(y)\subseteq \mathrm{im}(\varphi)$.
    Let $\gamma\in\mathrm{im}(\varphi)$, i.e.
    $$  \gamma = (\theta \circ g)(\eta,\sigma) \quad \text{for} \,\,\,\eta \in\Omega F, \,\,\, \sigma\in\Omega B.     $$
    Explicitly, we have
    $$   \varphi(\gamma)(s) = \begin{cases}
        \Phi(i\circ \eta(2s),2s) & 0\leq s\leq \tfrac{1}{2} \\ 
        \Phi(s\circ \sigma(2s-1),2-2s) & \tfrac{1}{2}\leq s \leq 1 .
    \end{cases}      $$
    Assume that there is a $t_*\in (0,1)$ with $\gamma(t_*) = p_0$.
    If $t_*<\tfrac{1}{2}$, then this implies
    $$    \Phi(i\circ \eta(2t_*),2t_*) = p_0     $$
    and thus
    $$         \Phi(p_0,2t_*)=  i\circ \eta(2t_*) \in V_1\times \{0\}    $$
    which contradicts \eqref{eq_property_cross_of_fiber_and_base}.
    Similarly, if $t_*\geq \tfrac{1}{2}$ then we obtain
    $$       \Phi(p_0,2-2t_*) =   s\circ \sigma(2t_*-1)\in \{0\}\times V_2     $$
    which again contradicts \eqref{eq_property_cross_of_fiber_and_base}.
    Hence, such a $t_*$ cannot exist.
    Consequently, all homology classes in $\mathrm{H}_{\bullet}(\Omega E,\{p_0\})$ have intersection multiplicity $1$.
 \end{proof}
 
 \begin{theorem}\label{tristan}
    Let $R$ be a commutative unital ring and take homology with coefficients in $R$.   
    Let $p\colon E\to B$ be a smooth fiber bundle with a section, where $B$ and $E$ are closed manifolds. Assume that both $B$ and the fiber $F$ are positive dimensional. Assume further that $E$ is $R$-oriented.
     Then the based string topology coproduct on $\mathrm{H}_{\bullet}(\Omega E,\{p_0\})$ vanishes.
 \end{theorem}
 
\begin{corollary}\label{isolde}
    Let $R$ be a commutative unital ring and take homology with coefficients in $R$.
    If a closed $R$-oriented manifold $M$ has non-vanishing based coproduct, then it cannot be the total space of a fiber bundle $M\to B$ with a section $s\colon B\to M$, where both fiber and base are positive dimensional.
\end{corollary}
 
 \begin{example}
     Consider an odd-dimensional sphere $\mathbb{S}^n$ with $n\geq 3$.
     Then its unit tangent bundle $U\mathbb{S}^n$ is a fiber bundle over $\mathbb{S}^n$ with fiber $\mathbb{S}^{n-1}$.
     The unit tangent bundle $U\mathbb{S}^n$ has a section since $\mathbb{S}^n$ admits a nowhere vanishing vector field.
     Consequently
     $$   \Omega \mathbb{S}^{n-1}\times \Omega\mathbb{S}^n \simeq \Omega (U\mathbb{S}^n)       $$
     and from the above corollary it follows that the based string topology coproduct vanishes on $\mathrm{H}_{\bullet}(\Omega (U\mathbb{S}^n),\{p_0\})$.
     Note that in this case the based string topology coproduct on both the base and the fiber does not vanish.
     Moreover, we remark that if $n = 4k+1$ for some $k\in\mathbb{N}$, then $U\mathbb{S}^n$ is not homotopy equivalent to the product $\mathbb{S}^n\times \mathbb{S}^{n-1}$, see \cite{james:1954}.
     Hence, the triviality of the based coproduct can not be deduced from Theorem \ref{theorem_int_mult_based} in this case and thus Theorem \ref{thm_fiber_bundle_with_section} is needed here.
 \end{example}

\subsection{Inclusions of submanifolds}

In the last part of this section we consider the situation of a closed manifold $M$ with a closed submanifold $N$.
The embedding $i\colon N\to M$ induces an inclusion of loop spaces $\Lambda N\hookrightarrow \Lambda M$.
In case that we consider the based loop space we choose the basepoint of $M$ such that it lies in $N$ and then also choose this point as a basepoint of $N$.
This yields an inclusion $\Omega N\hookrightarrow \Omega M$.

Equip $M$ with a Riemannian metric and let $E\to N$ be the normal bundle.
Then it is well-known that there is a tubular neighborhood $U$ of $N$, i.e. there is a homeomorphism $F\colon E\xrightarrow[]{} U$ such that the zero-section of $E$ is mapped to $N$.
Assume that there is a nowhere-vanishing section $\sigma \colon N\to E$ of the normal bundle.
Let $\Phi\colon N\times I\to M$ be the map
\begin{equation}\label{eq_phi}
      \Phi(p,s) = F ( s \cdot \sigma(p)) \quad \text{for}\,\,\,p\in N,\,\,s\in I .     
\end{equation}
 Then we have $\Phi(p,0) = p$ and $\Phi(p,s)\not\in N$ for $s> 0$.

 \begin{theorem}\label{theorem_int_mult_submanifold}
    Let $R$ be a commutative unital ring and take homology with coefficients in $R$. 
    Let $M$ be a closed manifold and $N$ a closed submanifold with normal bundle $E\to N$.
    Assume that the normal bundle has a section.
    Then for every homology class $X\in \mathrm{H}_{\bullet}(\Lambda N,N)$ it holds that
     $$    \mathrm{int}_{\Lambda M}(i_* X) = 1    $$
     where $i\colon \Lambda N\to \Lambda M$ is the inclusion.
 \end{theorem}
 \begin{proof}
     Let $\Phi\colon N\times I \to M$ be as in equation \eqref{eq_phi}.
     Let $X\in\mathrm{H}_{\bullet}(\Lambda N)$ be a homology class and $Y = i_* X$.
    Choose a continuous function $\rho\colon [0,1]\to [0,1]$ with $\rho(0) = \rho(1) = 0$ and $\rho(s)> 0$ for $s\in(0,1)$. 
     We define $\psi\colon \Lambda N\to \Lambda M$ by setting
     $$  \psi(\gamma)(t) = \Phi(\gamma(t),\rho(t) )   .  $$
     It is clear that $\psi\simeq i = \Phi(\cdot, 0)$ via the homotopy
     $$    H\colon \Lambda N\times I\to \Lambda M \quad H(\gamma,s)(t) = \Phi(\gamma(t),s\cdot\rho(t)) .  $$
     Therefore the classes $i_* X$ and $\psi_* X$ agree.
     Let $x$ be a representing cycle of the class $X$.
     Let $\gamma\in \mathrm{im}(\psi_* x)$, i.e. $\gamma= \psi(\sigma)$ for some $\sigma\in \Lambda N$.
     Assume that $\gamma$ has a non-trivial basepoint intersection, i.e. there is an $t_*\in (0,1)$ with $\gamma(0) = \gamma(t_*)$.
     But we have
     $$  \gamma(0) = \psi(\sigma)(0) \in N \quad \text{and} \quad \gamma(t_*) = \psi(\sigma)(t_*) \not\in N     $$
     by construction of $\psi$. This is a contradiction.
     Hence, all loops in $\mathrm{im}(\psi_* x)$ have only trivial self-intersections and therefore
     $$    \mathrm{int}_{\Lambda M}([\psi_* x])= \mathrm{int}_{\Lambda M}(Y) =    1  .  $$
     This shows the claim.
 \end{proof}
\begin{corollary}\label{cor_submanifold_coproduct}
    Let $R$ be a commutative unital ring and take homology with coefficients in $R$. 
    Let $M$ be a closed $R$-oriented manifold and $N$ a closed submanifold of $M$ such that the normal bundle has a nowhere vanishing section.
    Then the string topology coproduct vanishes on the subspace $U \subseteq \mathrm{H}_{\bullet}(\Lambda M,M)$ which is the image of the map induced by the inclusion
    $$  i_N\colon  \mathrm{H}_{\bullet}(\Lambda N) \xrightarrow[]{i_*} \mathrm{H}_{\bullet}(\Lambda M)\to \mathrm{H}_{\bullet}(\Lambda M,M) .    $$
\end{corollary}

\begin{proposition}
    Let $R$ be a commutative unital ring and take homology with coefficients in $R$. 
    Let $M_1,M_2$ be closed $R$-oriented manifolds with basepoints $p_0\in M_1$ and $q_0\in M_2$.
    Let 
     $j_1\colon M_1\hookrightarrow M = M_1\times M_2$ and $j_2\colon M_2\hookrightarrow M$ be the inclusions
     $$   j_1(r) = (r,q_0) \quad \text{and}\quad j_2(u) = (p_0,u)      $$
     and $i_1\colon \Lambda M_1\to \Lambda M$, resp. $i_2\colon \Lambda M_2\to \Lambda M$ be the induced maps between the loop spaces.
     Then all classes of the form
     $$   W = (i_1)_* X \quad \text{and} \quad Y = (i_2)_* Z     $$
     with $X\in\mathrm{H}_{\bullet}(\Lambda M_1)$ and $Z\in\mathrm{H}_{\bullet}(\Lambda M_2)$ have intersection multiplicity $1$ and thus trivial string topology coproduct.
\end{proposition}
\begin{proof}
    We only consider the case of $M_1\hookrightarrow M$.
    The other case is analogous.
    The inclusion $j_1$ is the inclusion of a submanifold and one checks that this has trivial normal bundle.
    Therefore we apply Theorem \ref{theorem_int_mult_submanifold}.
\end{proof}

Consider the free loop space of complex projective space $\mathbb{C}P^n$.
We have
    $$    \mathrm{H}_i(\Lambda\mathbb{C}P^k;\mathbb{Q}) \cong \mathbb{Q} \quad \text{for all}\,\,\, i\in\mathbb{N}_0,\,\,k\geq 1,   $$
    see e.g. \cite{stegemeyer:2022}.
\begin{theorem}\label{thm_cpn}
    Consider complex projective space $M = \mathbb{C}P^n$ and take $N = \mathbb{C}P^m$ as a submanifold which is embedded as the complex lines in $\mathbb{C}^{m+1}\subseteq \mathbb{C}^{n+1}$.
    Assume that $m<\tfrac{n}{2}$ and take homology with rational coefficients.
    Then the map
     $$    j_*\colon \mathrm{H}_{k}(\Lambda\mathbb{C}P^m) \to \mathrm{H}_{k}(\Lambda \mathbb{C}P^n)      $$
   induced by the inclusion $j\colon \Lambda \mathbb{C}P^m\hookrightarrow\Lambda \mathbb{C}P^m$ is trivial for degrees $k \not\in \{0,1,\ldots,2n\}\cup \{2n+2,2n+4,\ldots,4n\}$.
\end{theorem}
\begin{proof}
    Since $m<\tfrac{n}{2}$ the normal bundle of $N$ has rank $2n-2m > 2m$.
    Therefore the normal bundle of $N$ has a section, see e.g. \cite[Corollary VII.14.2]{bredon:2013}.
   If we compose the map 
   $$    j_*\colon \mathrm{H}_{\bullet}(\Lambda\mathbb{C}P^m) \to \mathrm{H}_{\bullet}(\Lambda \mathbb{C}P^n)   $$
  with the map induced by the inclusion of pairs $\iota\colon\mathrm{H}_{\bullet}(\Lambda \mathbb{C}P^n)\to \mathrm{H}_{\bullet}(\Lambda \mathbb{C}P^n,\mathbb{C}P^n)$ then we are in the situation of Corollary \ref{cor_submanifold_coproduct}.
    In particular the morphism $\iota$ is surjective and thus an isomorphism in degrees $i\not\in\{0,2,\ldots,2n\}$.
    By \cite{stegemeyer:2022} we know that the string topology coproduct is non-trivial for all non-trivial classes in $\mathrm{H}_i(\Lambda \mathbb{C}P^n,\mathbb{C}P^n)$ for $i\not\in I$ with 
    $$  I =  \{0,1,2,\ldots,2n\}\cup \{2n+2,2n+4,\ldots, 4n\} .   $$
    Let $i\in \mathbb{N}\setminus I$ and assume that $j_*\colon \mathrm{H}_i(\Lambda \mathbb{C}P^m) \to \mathrm{H}_i(\Lambda \mathbb{C}P^n)$ is non-zero.
    Then a generator $[X]\in\mathrm{H}_i(\Lambda \mathbb{C}P^m) $ has non-trivial image under the map $\iota\circ j_*$ in $\mathrm{H}_i(\Lambda\mathbb{C}P^n,\mathbb{C}P^n)$ and therefore non-trivial coproduct.
    This is a contradiction to Theorem \ref{theorem_int_mult_submanifold} and therefore we conclude that $j_*$ is trivial in all degrees $i\in \mathbb{N}\setminus I$.
\end{proof}

We want to conclude with a version of Theorem \ref{theorem_int_mult_submanifold} for the based loop space.
Note that we are not demanding that the tubular neighborhood of $N$ has a nowhere vanishing section.
\begin{theorem}\label{theorem_triviality_submfld_omega}
    Let $R$ be a commutative unital ring and take homology with coefficients in $R$.
    Let $M$ be a closed manifold and $N$ a closed submanifold such that $M$ and $N$ have the same basepoint.
    \begin{enumerate}
        \item For every homology class $X\in \mathrm{H}_{\bullet}(\Omega N,\{p_0\})$ it holds that
     $$    \mathrm{int}_{\Omega}(i_* X) = 1    $$
     where $i\colon \Omega N\to \Omega M$ is the inclusion.
     \item Assuming that $M$ is $R$-oriented the based string topology coproduct vanishes on the subspace $U \subseteq \mathrm{H}_{\bullet}(\Omega M,\{p_0\})$ which is the image of the map induced by the inclusion
    $$    \mathrm{H}_{\bullet}(\Omega N) \xrightarrow[]{i_*} \mathrm{H}_{\bullet}(\Omega M)\to \mathrm{H}_{\bullet}(\Omega M,\{p_0\}) .    $$
    \end{enumerate}
\end{theorem}
\begin{proof}
    Let $k = \mathrm{dim}(N)$ and $n = \mathrm{dim}(M)$.
    Let $p_0\in N$ be the basepoint in $N$ which is mapped to the basepoint $r_0\in M$ under the inclusion $N\hookrightarrow M$.
    Let $U\subseteq M$ be a neighborhood of $r_0$ such that there is a homeomorphism $\varphi\colon U\to B_1^n(0)$, where $B_1^n(0)$ is the ball with radius $1$ in $\mathbb{R}^n$.
    Moreover, we can choose $U$ and $\varphi$ such that
    $$  \varphi(N\cap U) =  B_1^n(0) \cap  \mathbb{R}^k .     $$
    In particular this induces a map $\varphi'\colon N\cap U\to B_1^k(0)$.
    Let $\sigma\colon [0,1)\to [0,1]$ be a smooth function which is equal to $1$ on $[0,\tfrac{1}{2}]$ and vanishes for $|t|\geq \tfrac{3}{4}$.
    Then we define a map
    $$   \Phi\colon N\times I \to M    $$
    by setting
    $$    \Phi(p,s) = \begin{cases}
        \varphi^{-1}(\varphi'(p),s\cdot \sigma(||\varphi'(p)||),0,\ldots,0)  &  p\in N\cap U \\ 
        p & \text{else} .        
    \end{cases}      $$
    Note that we put $s\cdot \sigma(||\varphi'(p)||)$ in the $(k+1)$-st component in the above equation.
    Choose a smooth function $\rho\colon [0,1]\to[0,1]$ with $\rho(0)=\rho(1)= 0$ and $\rho(s) > 0$ for $s\in(0,1)$.
    We define $\psi\colon \Omega N\to \Omega M$ by
    $$    \psi(\gamma)(t) =  \Phi(\gamma(t),\rho(t))  .   $$
    We see that $\psi\simeq i$ via the homotopy
    $$   H\colon \Omega N\times I\to \Omega M,\quad H(\gamma,s)(t) = \Phi(\gamma(t),s\cdot \rho(t)) .    $$
    Hence, if $X\in \mathrm{H}_{\bullet}(\Omega N)$ is a homology class then $i_* X = \psi_* X$.
    With the same consideration as in the proof of Theorem \ref{theorem_int_mult_submanifold} we see that
    $$   \mathrm{int}_{\Omega M}\left(\iota\circ i_* (X)\right) =  1 ,  $$
    where $\iota$ is the projection $\mathrm{H}_{\bullet}(\Omega M)\to \mathrm{H}_{\bullet}(\Omega M,\{p_0\})$. The second part of the statement then follows easily.
\end{proof}

\bibliography{kupper_stegemeyer_references}  
\bibliographystyle{spmpsci}

\end{document}